\documentclass[11pt]{article}%
\usepackage{amsmath}
\usepackage{amsfonts}
\usepackage{amssymb}
\usepackage{graphicx}
\usepackage{a4wide}

\providecommand{\U}[1]{\protect\rule{.1in}{.1in}}
\newtheorem{theorem}{Theorem}

\newtheorem{example}[theorem]{Example}

\begin{document}

\title{On a variant of Schu's lemma}
\author{{R. I. Bo\c{t}\thanks{Faculty of Mathematics, University of Vienna, Vienna,
Austria, email: \texttt{radu.bot@univie.ac.at.}}} ~and C.
Z\u{a}linescu\thanks{Octav Mayer Institute of Mathematics, Ia\c{s}i Branch of
Romanian Academy, Ia\c{s}i, Romania, email: \texttt{zalinesc@uaic.ro.}}}
\date{}
\maketitle

\begin{abstract}
In this note, we demonstrate that an incorrect statement has been propagated
in multiple papers, stemming from the substitution of ``lim'' with ``limsup''
for a sequence in Lemma 1.3 of the paper [J. Schu: Weak and strong convergence
to fixed points of asymptotically nonexpansive mappings,
Bull.\ Austral.\ Math.\ Soc.\ 43 (1991), 153--159]. This occurred over a span
of more than 20 years, with the earliest paper we identified using this
incorrect statement dating back to 2002.

\end{abstract}

In Schu's paper \cite{Sch91}, the following result is stated without
proof, {with ``[8]'' referring to Zeidler's book, which we cite as
\cite{Zei86}:}

\smallskip``LEMMA 1.3. (compare [8, p.\ 484]) \textit{Let $(E,\left\Vert
\cdot\right\Vert )$ be a uniformly convex Banach space, $0<b<c<1$, $a\geq0$,
$(t_{n})\subset\lbrack b,c]^{\mathbb{N}}$ and $(x_{n}),(y_{n})\in
E^{\mathbb{N}}$ such that $\limsup\left\Vert x_{n}\right\Vert \leq a$,
$\limsup\left\Vert y_{n}\right\Vert \leq a$ and $\lim\left\Vert t_{n}%
x_{n}+(1-t_{n})y_{n}\right\Vert =a$. Then $\lim\left\Vert x_{n}-y_{n}%
\right\Vert =0$}."

\medskip
{Interestingly, Dotson had presented a similar result 21 years
earlier, but most of the literature we know of using this result
attributes it to Schu. On page 68 of \cite{Dot70}, Dotson states the
following:}

\smallskip

\textquotedblleft The following lemma is an easy consequence of uniform convexity.

\smallskip LEMMA 3. \textit{Suppose $E$ is a uniformly convex Banach space.
Suppose $0<a<b<1$, and $\{t_{n}\}$ is a sequence in $[a,b]$. Suppose
$\{w_{n}\}$, $\{y_{n}\}$ are sequences in $E$ such that $\left\Vert
w_{n}\right\Vert \leq1$, $\left\Vert y_{n}\right\Vert \leq1$ for all $n$.
Define $\{z_{n}\}$ in $E$ by $z_{n}=(1-t_{n})w_{n}+t_{n}y_{n}$. If
$\lim\left\Vert z_{n}\right\Vert =1$, then $\lim\left\Vert w_{n}%
-y_{n}\right\Vert =0$.}"

\smallskip Note that \cite[Lemma 3]{Dot70} is recalled in \cite[page
376]{SenDot74}.

\medskip Taking $t_{n}:=t\in(0,1)$ for $n\in\mathbb{N}^{\ast}$
$(:=\{1,2,...,n,...\})$, \cite[Lemma 1.3]{Sch91} reduces to {Problem
10.1 (c)} from \cite[p.\ 484]{Zei86}. {Additionally, note that} the
version of \cite[Lemma 1.3]{Sch91} {where the interval $[b,c]$ is
replaced by $[\varepsilon,1-\varepsilon]$ with
$\varepsilon\in(0,1)$, } is attributed also to
\textquotedblleft(Zeidler (1986))" in \cite[Lemma 2.1]{TakRed19}.

{In \cite{KhAhAb23}, a variant of Schu's lemma, formulated with a
\textbf{weaker hypothesis}, is introduced as follows:}

\smallskip\textquotedblleft In 1991, Sahu [21] established an important
property of UCBS, which can be stated as follows.

\smallskip Lemma 2.2. [21] \textit{Assume that $\Omega$ is a UCBS and
$\{z_{\ell}\}$ be a sequence in $(0,1)$ for all $\ell\geq1$. If $\{a_{\ell}\}$
and $\{b_{\ell}\}$ are in $\Omega$ such that $\limsup_{\ell\rightarrow\infty
}\left\Vert a_{\ell}\right\Vert \leq m$, $\limsup_{\ell\rightarrow\infty
}\left\Vert b_{\ell}\right\Vert \leq m$ and $\lim_{\ell\rightarrow\infty
}\left\Vert (1-z_{\ell})a_{\ell}+z_{\ell}b_{\ell}\right\Vert =m$ for some
$m\geq0$. Then $\lim_{\ell\rightarrow\infty}\left\Vert a_{\ell}-b_{\ell
}\right\Vert =0$.}"

\smallskip In fact, reference [21] from the preceding quoted text is
{none other than Schu's paper} \cite{Sch91}. {Additionally, note
that} Lemma 2.5 from \cite{ReIqBaEtBo23}, also attributed to
\cite{Sch91}, is equivalent to \cite[Lemma 2.2]{KhAhAb23}.

{Note that this variant of Schu's lemma \cite[Lemma 1.3]{Sch91} is
\textbf{false}, as illustrated by the following example.}

\begin{example}
\label{ex0}Let $(X,\left\Vert \cdot\right\Vert )$ be a non-trivial normed
vector space, $x\in X\setminus\{0\}$, $x_{n}:=x$ and $y_{n}:=-x$ for
$n\in\mathbb{N}^{\ast};$ moreover, consider $\{t_{n}\}\subset(0,1)$ with
$t_{n}\rightarrow0$. Then $\left\Vert x_{n}\right\Vert =\left\Vert
y_{n}\right\Vert =\left\Vert x\right\Vert $, and so $\limsup_{n\rightarrow
\infty}\left\Vert x_{n}\right\Vert =\limsup_{n\rightarrow\infty}\left\Vert
y_{n}\right\Vert =\left\Vert x\right\Vert >0$; moreover, $t_{n}x_{n}%
+(1-t_{n})y_{n}=(2t_{n}-1)x$, and so $\lim_{n\rightarrow\infty}\left\Vert
t_{n}x_{n}+(1-t_{n})y_{n}\right\Vert =\lim_{n\rightarrow\infty}\left\vert
2t_{n}-1\right\vert \cdot\left\Vert x\right\Vert =\left\Vert x\right\Vert .$
It is obvious that $\lim_{n\rightarrow\infty}\left\Vert x_{n}-y_{n}\right\Vert
=2\left\Vert x\right\Vert \neq0$.
\end{example}

In \cite{LaoPan10}, {Laowang and Panyanak introduced an extension of
Schu's lemma to the setting of uniformly convex hyperbolic spaces,
first announced in \cite[Lemma 2.7]{LaoPan09} for CAT$(0)$ spaces,
as follows:}

\smallskip\textquotedblleft The following result is a characterization of
uniformly convex hyperbolic spaces which is an analog of Lemma 1.3 of Schu
[25]. It can be applied to a CAT$(0)$ space as well.

\smallskip

Lemma 2.9. \textit{Let $(X,d,W)$ be a uniformly convex hyperbolic space with
modulus of convexity $\eta$, and let $x\in X$. Suppose that $\eta$ increases
with $r$ (for a fixed $\varepsilon$) and suppose that $\{t_{n}\}$ is a
sequence in $[b,c]$ for some $b,c\in(0,1)$, and $\{x_{n}\}$, $\{y_{n}\}$ are
sequences in $X$ such that $\limsup_{n\rightarrow\infty}d(x_{n},x)\leq r$,
$\limsup_{n\rightarrow\infty}d(y_{n},x)\leq r$, and $\lim_{n\rightarrow\infty
}d((1-t_{n})x_{n})\oplus t_{n}y_{n},x)=r$ for some $r\geq0$. Then
$\lim_{n\rightarrow\infty}d(x_{n},y_{n})=0$. (2.17)}"

\medskip

The proof of \cite[Lemma 2.9]{LaoPan10} is an adaptation to this context of
that for Problem 10.1 (c) in Zeidler's book \cite{Zei86}.

The formulation and the proof of \cite[Lemma 2.9]{LaoPan10} are replicated in
\cite[Lemma 4.5]{NanPan10}, while, for $t_{n}:=1/2$ for $n\in\mathbb{N}^{\ast
}$, it can be found as Lemma 2.2 in \cite{KhaKha11}.

Lemma 2.9 from \cite{LaoPan10} is (practically) replicated also in \cite[Lemma
2.5]{KhFuKh12}, where an \emph{alternative proof} is provided. However, this
proof contains the following \emph{false assertion}: \textquotedblleft Since
{$\limsup_{n\rightarrow\infty}d(x_{n},x)\leq r$ and $\limsup_{n\rightarrow
\infty}d(y_{n},x)\leq r$} we have: (i)~{$d(x_{n},x)\leq r+\frac{1}{n}$;
(ii)~$d(y_{n},x)\leq r+\frac{1}{n}$ }for each $n\geq1$".

Lemma 2.9 from \cite{LaoPan10} is also presented as Lemma 2.5 in
\cite{FukKha14}, with a reference to \cite[Lemma 2.3]{FukKha14}, and as Lemma
3.2 in \cite{ChWaWaTaMa14}.

Note that \cite[Lemma 2.9]{LaoPan10} is stated and proved for $X$ as a
CAT$(0)$ space in \cite[Lemma 3.2]{ChWaLeChYa12} under the supplementary
condition $0<b(1-c)\leq\tfrac{1}{2}$. This condition appears in many papers
from our bibliography. It is easy to observe that this condition is
superfluous. Indeed, $b,c\in(0,1)$ implies that $0<b(1-c)$, while the fact
that (there exists) $t_{n}\in\lbrack b,c]$ implies that $b\leq c$; hence
$0<b\leq c<1$. Because $b^{2}+(b-1)^{2}=2b^{2}-2b+1 > 0$ and $b>0$, it follows
that $\frac{1}{2b}>1- b \geq1-c $. Consequently, $b(1-c)<\frac{1}{2}$.

\medskip

Furthermore, an extension of Schu's lemma to the setting of modular
function spaces was introduced in \cite{DehKoz12} (see also
\cite[Lemma 4.2]{KhaKoz15}):

\smallskip\textquotedblleft Lemma 3.2. \textit{Let $\rho\in\Re$ be (UUC1) and
let $\{t_{k}\}\subset(0,1)$ be bounded away from $0$ and $1$. If there exists
$R>0$ such that $\limsup_{n\rightarrow\infty}\rho(f_{n})\leq R$,
$\limsup_{n\rightarrow\infty}\rho(g_{n})\leq R$, $\lim_{n\rightarrow\infty
}\rho(t_{n}f_{n}+(1-t_{n})g_{n})=R$,\ then $\lim_{n\rightarrow\infty}%
\rho(f_{n}-g_{n})=0.$}"

\medskip

\medskip In the following we present an \emph{alternative proof for
\cite[Lemma 1.3]{Sch91}}:

Since $\limsup\left\Vert x_{n}\right\Vert \leq a$ and $\limsup\left\Vert
y_{n}\right\Vert \leq a$, there exists $R>0$ such that $\left\Vert
x_{n}\right\Vert ,\left\Vert y_{n}\right\Vert \leq R$ for every $n\in
\mathbb{N}^{*}$. $(E,\left\Vert \cdot\right\Vert )$ being uniformly convex, by
\cite[Theorem 4.1 (ii)]{Zal02}, $\left\Vert \cdot\right\Vert ^{2}$ is a
uniformly convex function on $B:=\{x\in E\mid\left\Vert x\right\Vert \leq
R\}$. This means that there exists a strictly increasing function
$\rho:\mathbb{R}_{+} \rightarrow\mathbb{R}_{+}$, with $\rho(0)=0$, such that
\[
\left\Vert \lambda x+(1-\lambda)y\right\Vert ^{2}\leq\lambda\left\Vert
x\right\Vert ^{2}+(1-\lambda)\left\Vert y\right\Vert ^{2}-\lambda
(1-\lambda)\rho\left(  \left\Vert x-y\right\Vert \right)  \quad\forall x,y\in
B ,\ \forall\lambda\in\lbrack0,1].
\]
Setting $\gamma:=\min\{b,1-c\}$ $(>0)$, one has $\lambda(1-\lambda)\geq
\gamma^{2}$, and so
\begin{equation}
\gamma^{2}\rho\left(  \left\Vert x_{n}-y_{n}\right\Vert \right)  \leq
t_{n}\left\Vert x_{n}\right\Vert ^{2}+(1-t_{n})\left\Vert y_{n}\right\Vert
^{2}-\left\Vert t_{n}x_{n}+(1-t_{n})y_{n}\right\Vert ^{2}\quad\forall
n\in\mathbb{N}. \label{r1}%
\end{equation}
Assume that the conclusion does not hold. Then there exist a
strictly increasing sequence $(n_{k})\subset\mathbb{N}^{\ast}$ and
$\varepsilon>0$ such that $t_{n_{k}}\rightarrow t\in\lbrack b,c]$ as
$k \to\infty$, and $\left\Vert x_{n_{k} }-y_{n_{k}}\right\Vert
\geq\varepsilon$ for every $k\in \mathbb{N}^{\ast}$. {Replacing, if
necessary, $n$ by $n_{k}$ in (\ref{r1}), using that $\rho\left(
\left\Vert x_{n_{k}}-y_{n_{k}}\right\Vert \right)  \geq
\rho(\varepsilon)$ for every $k\in\mathbb{N}^{\ast}$, and passing to
$\limsup$
in (\ref{r1}) as $k \to\infty$,} one gets the contradiction $\gamma^{2}%
\rho(\varepsilon)\leq ta^{2}+(1-t)a^{2}-a^{2}=0$. \hfill$\square$

\medskip

{In \cite{KiKiTa02}, Kim, Kiuchi, and Takahashi introduced a variant
of Schu's lemma in uniformly convex Banach spaces, assuming a
\textbf{weaker hypothesis}, with ``[9]'' referring to \cite{Sch91}:}

\smallskip

\textquotedblleft Lemma 2.4 ([9]). \textit{Let $E$ be a uniformly convex
Banach space, let $0<b\leq t_{n}\leq c<1$ for all $n\in N$, and let
$\{x_{n}\}$ and $\{y_{n}\}$ be sequences of E such that $\overline{\lim
}\left\Vert x_{n}\right\Vert \leq a$, $\overline{\lim}\left\Vert
y_{n}\right\Vert \leq a$, and $\overline{\lim}\left\Vert t_{n}x_{n}%
+(1-t_{n})y_{n}\right\Vert =a$ for some $a\geq0$. Then, it holds that
$\lim_{n\rightarrow\infty}\left\Vert x_{n}-y_{n}\right\Vert =0$.}"

\smallskip

{No proof is given for \cite[Lemma 2.4]{KiKiTa02}. However, Agarwal,
O'Regan, and Sahu formulated the following result in their book
\cite{AgOrSa09}, for which a proof is also given:}

\smallskip

\textquotedblleft Theorem 2.3.13 \textit{Let $X$ be a uniformly convex Banach
space and let {$\{t_{n}\}$} be a sequence of real numbers in $(0,1)$ bounded
away from $0$ and $1$. Let $\{x_{n}\}$ and {$\{y_{n}\}$} be two sequences in
$X$ such that
\[
{\limsup_{n\rightarrow\infty}\left\Vert x_{n}\right\Vert \leq a,\ \ {\limsup
_{n\rightarrow\infty}}\left\Vert y_{n}\right\Vert \leq a}\text{\ \ {and}%
}{\ ~\limsup_{n\rightarrow\infty}\left\Vert t_{n}x_{n}+(1-t_{n})y_{n}%
\right\Vert =a}%
\]
for some $a\geq0$. Then $\lim\left\Vert x_{n}-y_{n}\right\Vert =0$.}"

\medskip

{Replacing $a$ by $r$ in \cite[Theorem 2.3.13]{AgOrSa09} (above) and
the pair $(\alpha,\beta)$ from its proof by $(p,q)$ one gets the
statement of \cite[Theorem 2.7]{Des20}; compare $R\in(a,a+1)$ in its
proof with $r\in(a,a+1)$ from the proof of \cite[Theorem
2.3.13]{AgOrSa09}. Notice that \cite{AgOrSa09} is mentioned in the
bibliography of \cite{Des20}, but is not cited in this context.}

\smallskip In fact \cite[Lemma 2.4]{KiKiTa02} (as well as \cite[Theorem 2.3.13]%
{AgOrSa09}) is \textbf{false}, as demonstrated by the following
example.

\begin{example}
\label{ex1}Let $(X,\left\Vert \cdot\right\Vert )$ be a non-trivial normed
vector space, $x\in X$ with $\left\Vert x\right\Vert =1$ $(=:a)$, $x_{n}:=x,$
$y_{2n}:=0,$ $y_{2n-1}:=x,$ $t_{n}:=\lambda\in(0,1)$ for $n\in\mathbb{N}%
^{\ast}$. Then $\limsup_{n\rightarrow\infty}\left\Vert x_{n}\right\Vert =1,$
$\limsup_{n\rightarrow\infty}\left\Vert y_{n}\right\Vert =1,$ $t_{n}%
x_{n}+(1-t_{n})y_{n}=\lambda x$ for even $n$, $t_{n}x_{n}+(1-t_{n})y_{n}=x$
for odd $n$, and so $\limsup_{n\rightarrow\infty}\left\Vert t_{n}%
x_{n}+(1-t_{n})y_{n}\right\Vert =\limsup_{n\rightarrow\infty}1=1.$ It is
obvious that $\lim_{n\rightarrow\infty}\left\Vert x_{n}-y_{n}\right\Vert $
does not exist, and so it is different from $0$.
\end{example}

{Below is a list of statements equivalent to \cite[Lemma
2.4]{KiKiTa02}, given obviously without proofs, and used in the
proofs of other results in the corresponding works. In most cases,
\cite[Lemma 1.3]{Sch91} is cited as the reference for this result;
exceptions will be noted explicitly. For the interval $[b,c]$ for
some $0<b\leq c<1$ often the supplementary condition(s) $b<c$ and/or
$0<b(1-c)\leq\tfrac{1}{2}$ are required, or this is replaced by
$[\varepsilon,1-\varepsilon]$ with $\varepsilon\in(0,1)$:}

\medskip\noindent\cite[Lemma 2.3]{KiKiTa04}, \cite[Lemma 2.1]{Sha05},
\cite[Lemma 2.3]{JeoKim06}, \cite[Lemma 2.4]{ShaAld06}, \cite[Lemma
2.1]{ShaUdo06}, \cite[Lemma 5]{ChiOfo07}, \cite[Lemma
2.1]{KiOzAk07}, \cite[Lemma 2.2]{ZhoWan07}, \cite[Lemma
2.1]{DenLiu08} (attributed to \cite{Xu91}), \cite[Lemma 2.1]{Mou08},
\cite[Lemma 2.2]{Yao08}, \cite[Lemma 2.2]{YaoWan08}, \cite[Lemma
2.2]{YaoZho08}, \cite[Lemma 2.2]{PenYao09}, \cite[Lemma
2.3]{PenYao09b}, {\cite[Theorem 4.3.1]{Riz09}}, \cite[Lemma
2.2]{Tem09}, \cite[Lemma 2.3]{Tem09b}, \cite[Lemma 2.2]{YaoYan09}
(attributed to \cite{Sch91b}), \cite[Lemma 2.6]{PenYao10},
\cite[Lemma 2.2]{Tem10}, \cite[Lemma 2.1]{ManTak11}, \cite[Lemma
2.1]{AkSoCh12}, \cite[Lemma 2.4]{Jun12}, \cite[Lemma 2.2]{MaChe12},
\cite[Lemma 2.7]{Sal12}, \cite[Lemma 2]{KarOzd13}, \cite[Lemma
2.3]{Tem13}, \cite[Lemma 2.2]{Tem13b}, \cite[Lemma 6]{WaLiKo13},
\cite[Lemma 1.3]{HouDu14}, \cite[Lemma 2.9]{JhaSal14}, \cite[Lemma
3.2]{ThThPo14}, \cite[Lemma 1.9]{JamAbe15} (\cite{AgOrSa09}),
\cite[Lemma 2.2]{Jun15}, \cite[Lemma 1.4]{RazMor15}, \cite[Lemma
2.2]{Tem15},
\cite[Lemma 12]{Abd16}, \cite[Lemma 4.2]{AsKhAl16}, \cite[Lemma 2.1]%
{SaPaTi16}, \cite[Lemma 2.2]{Shr16}, \cite[Lemma 2.4]{ThThPo16}, \cite[Lemma
2.4]{ThThPo16b}, \cite[Lemma 2.10]{Tri16}, \cite[Lemma 2.2]{WoGoHa16}
(attributed to \cite{Zei86}), \cite[Lemma 1.4]{PanBho17}, \cite[Lemma
2.3]{UllArs17}, \cite[Lemma 2.2]{AbdRas18}, \cite[Lemma 2.7]{AlKoTr18},
\cite[Lemma 2.4]{HuUlAr18}, \cite[Lemma 2]{SuChSu18}, \cite[Lemma
2.4]{UllArs18b}, \cite[Lemma 2.4]{UllArs18}, \cite[Lemma 2.1]{WatKla18},
\cite[Lemma 2.6]{AliAli19}, \cite[Lemma 3]{AlAlKu19}, \cite[Lemma
1.4]{BhuTiw19}, \cite[Lemma 2.5]{FeShCh19}, \cite[Lemma 1.7]{PanBho19} (with
$0<g_{n}<1$), \cite[Lemma 2.1]{SrSuSuCh19}, \cite[Lemma 1.4]{ThiYam19},
\cite[Lemma 1]{UsuPos19}, \cite[Lemma 4.6]{AlGuAk20}, \cite[Lemma
2.6]{AlAlNi20}, \cite[Lemma 5]{BejPos20}, \cite[Lemma 2.4]{ChYaTh20},
\cite[Lemma 3]{CioTur20}, \cite[Lemma 2.1]{GarUdi20}, \cite[Lemma
2.1]{GarUdi20b}, \cite[Lemma 2.2]{GaUdKh20}, \cite[Lemma 3]{HaDeAgAlHu20},
\cite[Lemma 3]{HouTur20}, \cite[Lemma 2.2]{Man20}, \cite[Lemma 2.8]{RanArt20},
\cite[Lemma 2]{RimKim20}, \cite[Lemma 2.2]{UsBePo20}, \cite[Lemma 3]%
{YaAuTh20}, \cite[Lemma 2]{AlAlKh21}, \cite[Lemma 2]{GaAbUd21}, \cite[Lemma
4]{HuHuAl21}, \cite[Lemma 2.2]{KhaUdd21}, \cite[Lemma 3]{KhUdAlGe21},
\cite[Lemma 4.2.4]{Kil21}, \cite[Lemma 2.16]{OfeIgb21}, \cite[Lemma
2.14]{OfUdIg21}, \cite[Lemma 3]{OfUdIg21b}, \cite[Lemma 1.2]{Oze21}
(\cite{SenDot74}), \cite[Lemma 2.10]{RaDiBa21}, \cite[Lemma 5]{ThKaPhSuIn21},
\cite[Lemma 2.7]{UdOfIg21}, \cite[Lemma 4]{AlHuKaCh22}, \cite[Lemma
2.4]{AlJuAl22}, \cite[Lemma 2.3]{AlGaUdNi22}, \cite[Lemma 1]{BeAbAs22} (with
$\rho_{n} \in(0,1)$), \cite[Lemma 2.7]{BejCio22}, \cite[Lemma 3]{BejPos22},
\cite[Lemma 2.11]{DewGur22}, \cite[Lemma 2.3]{GaUdBa22}, \cite[Lemma
2]{GoPrDe22} (attributed to \cite{ThThPo16}), \cite[Lemma 6]{JiShAhShBo22},
\cite[Lemma 12]{JuAlKu22}, \cite[Lemma 2]{MaGuAtAb22}, \cite[Lemma
1]{OfHuJoUdIsSuCh22}, \cite[Lemma 2.6]{OfIsAlAh22}, \cite[Lemma 2.6]%
{OkOfIs22}, \cite[Lemma 2.4]{OkeUgw22}, \cite[Lemma 2.2]{PaRoPaSaPa22},
\cite[Lemma 2.10]{SaBaGu22}, \cite[Lemma 1.4]{SalMai22} (with $t_{n}\in[0,1]$,
\cite{AgOrSa09}), \cite[Lemma 7]{ThInSuPh22}, \cite[Lemma 2.1]{UdGaAbMl22},
\cite[Lemma 2.2]{UsuPos22}, \cite[Lemma 3]{UsTuPo22}, \cite[Lemma 1]%
{YamThi22}, \cite[Lemma 1]{AlAlDaTuZaMa23}, \cite[Lemma
1]{AliJub23}, \cite[Lemma 1.6]{BhKuKu23}, \cite[Lemma
2.3]{ChYaTh23}, \cite[Theorem 2.6]{DeGoRa23} (\cite{Zei86}; {compare
with \cite{Des20}}), \cite[Lemma 2.5]{DewGur23}, \cite[Lemma
2]{DeGuAhAlMl23}, \cite[Lemma 2.3]{EkiTem23}, \cite[Lemma
2]{GauKau23}, \cite[Lemma 5]{HamKat23}, \cite[Lemma 5]{HamKat23b},
\cite[Lemma 1]{KhUdSw23}, \cite[Lemma 4]{KiYaChTh23}, \cite[Lemma
(1.2)]{MaiSal23} (\cite{AgOrSa09}), \cite[Lemma 2.1]{Sae23},
\cite[Lemma 2.12]{SahBan23}, \cite[Lemma 2.1]{SaBaGu23}, \cite[Lemma
1.3]{Tem23}, \cite[Lemma 1]{Tem23b}, \cite[Lemma 2.4]{TemZin23},
\cite[Lemma 2.3]{UlSaBiAhIbJa23}, \cite[Lemma 2.1]{AhShNa24},
\cite[Lemma 2.2]{AlHaReSe24}, \cite[Lemma 2.4]{Dew24}, \cite[Lemma
2.1]{DewGur24}, \cite[Lemma 2.2]{DeRaMiRa24}, \cite[Lemma
2.3]{GauKau24}, \cite[Lemma 1]{IqAlSuHu24}, \cite[Lemma
2.1]{MaBaAlAhAr24} (with $0<t_{s}<1$),
\cite[Lemma 1]{Nav24}, \cite[Lemma 1]{NaUlGd24}, \cite[Lemma 2.1]%
{OkUdAlKaAh24}, \cite[Lemma 3]{SahBan24}, \cite[Lemma 2.6]{SahBan24b},
\cite[Lemma 1]{UgObUdNaOfNa24}, \cite[Lemma 2.5]{UgObUdOkChNa24}, \cite[Lemma
3]{YamThi24}.

\medskip
{In the setting of CAT(0) spaces, Zhou and Cui formulated the
following result in \cite{ZhoCui15}, which assumes a \textbf{weaker
hypothesis} than that in \cite[Lemma 2.9]{LaoPan10}, recalled
above:}

\smallskip\textquotedblleft Lemma 2.2 ([6, Lemma 2.9]). \textit{Let $(M,d)$ be
a CAT$(0)$ space and let $x\in M$. Suppose {$\{t_{n}\}$} is a sequence in
$[b,c]$ for some $b,c\in(0,1)$ and {$\{x_{n}\},$ $\{y_{n}\}$} are sequences in
$M$ such that $\limsup_{n\rightarrow\infty}d(x_{n},x)\leq r$, $\limsup
_{n\rightarrow\infty}d(y_{n},x)\leq r$ and $\limsup_{n\rightarrow\infty
}d((1-t_{n})x_{n})\oplus t_{n}y_{n},x)=r$ for some $r\geq0$. Then
$\lim_{n\rightarrow\infty}d(x_{n},y_{n})=0$.}"

\medskip The reference ``[6]'' mentioned in \cite[Lemma 2.2]{ZhoCui15} is our
reference \cite{DhKiSi06} and does not contain a statement called ``Lemma
2.9''. Of course, no proof is given for this ``stronger'' version of
\cite[Lemma 2.9]{LaoPan10} in \cite{ZhoCui15}. Since any Hilbert space is
CAT$(0)$, Example \ref{ex1} shows that \cite[Lemma 2.2]{ZhoCui15} is
\textbf{false}.

\medskip Below is a list of statements equivalent to \cite[Lemma
2.2]{ZhoCui15} used in the proofs of other results in the corresponding works.
For each paper, we add the version of \cite[Lemma 2.9]{LaoPan10} that is
cited, or just the work if a particular result is not mentioned:

\medskip\noindent\cite[Lemma 2.3]{ThThPo15} (\cite[Lemma 4.5]{NanPan10}),
\cite[Lemma 2.2]{Zha15} (\cite{KhFuKh12,Leu10}), \cite[Lemma
2.5]{Ali16} (\cite[Lemma 2.2]{KhaKha11}), \cite[Lemma
1.12]{PanSin16} (\cite[Lemma 4.5]{NanPan10}), \cite[Lemma
1.11]{SuaKli16} (\cite[Lemma 2.5]{KhFuKh12} via \cite[Lemma
1.3]{ChWaWaTaMa14}), \cite[Lemma 2]{AbkShe17} (\cite[Lemma
4.5]{NanPan10}), \cite[Lemma 2.11]{QiaDen17} (\cite{LaoPan10}),
\cite[Lemma 2.5]{RasAbk17} (\cite[Lemma 4.5]{NanPan10}), \cite[Lemma
2.2]{AbkRas18} (\cite[Lemma 4.5]{NanPan10}), \cite[Lemma
2.7]{UlIqAr18} (\cite[Lemma 2.9]{LaoPan10}), \cite[Lemma
2.12]{UlKhAr18} (\cite{Noo00}), \cite[Lemma 2.6]{BaSuKh20}
(\cite[Lemma 2.5]{KhFuKh12} via \cite[Lemma 1.11]{SuaKli16}),
\cite[Lemma 2.11]{DaPaSa20} (\cite[Lemma 2.5]{KhFuKh12}),
\cite[Lemma 2.5]{KuPaKuCh20} (\cite{LaoPan09}), \cite[Lemma
5]{AbIqDeAh21} (\cite{LaoPan10}), \cite[Lemma
2.2]{AbRaGu21} (\cite{NanPan10}), \cite[Lemma 1.4]{KimKan21} (\cite{KhFuKh12}%
), \cite[Lemma 5.6]{LamPan21} with {$\limsup_{n\rightarrow\infty}d(t_{n}%
x_{n}\oplus(1-t_{n})y_{n})\leq r$} (\cite{LaoPan10}), \cite[Lemma 2.1]{Hao22}
(\cite[Lemma 2.5]{KhFuKh12}), \cite[Lemma 3.8]{KiDaPaSa22} (\cite{KhFuKh12}),
\cite[Lemma 2.14]{SaKuKa21} (\cite{DhKiSi06}), \cite[Lemma 3]{TufZeg22}
(\cite{LaoPan10}), \cite[Lemma 2.10]{AbAhSh23} (\cite{LaoPan10}), \cite[Lemma
12]{RaArThKe23} (\cite{Sch91}), \cite[Lemma 2.10]{DaTiSh24} (\cite{LaoPan09}),
\cite[Lemma 2.3]{KiDaPaVe24} with {$\limsup_{n\rightarrow\infty}%
d(W(x_{n},y_{n},t_{n}),x)\leq c$} (\cite{KhFuKh12}), \cite[Lemma
2.8]{ShrCha24} (\cite{KhFuKh12}), \cite[Lemma 5]{TaAhZaRaIsNa24} (\cite{Sch91}).

\medskip
{As for the setting of modular function spaces, Bejenaru and
Postolache formulated the following result in \cite{BejPos19}, which
assumes a \textbf{weaker hypothesis} than that in \cite[Lemma
3.2]{DehKoz12}, recalled above:}

\smallskip\textquotedblleft Lemma 1. \textit{Suppose that $\rho$ satisfies
property (UUC1) and let $\{\alpha_{l}\}\subset\lbrack a,b]$, where $0<a<b<1$.
If there exists a positive real number $r$ such that $\limsup_{l\rightarrow
\infty}\rho(\alpha_{l}x_{l}+(1-\alpha_{l})y_{l})=r$, $\limsup_{l\rightarrow
\infty}\rho(x_{l})\leq r$ and $\limsup_{l\rightarrow\infty}\rho(y_{l})\leq r$,
then $\lim_{l\rightarrow\infty}\rho(x_{l}-y_{l})=0$ ([8], cf.\ [5]).}"

\medskip The references ``[8]'' and ``[5]'' mentioned in \cite[Lemma
1]{BejPos19} are our references \cite{Kha18} and \cite{DehKoz12},
respectively, where \cite[Lemma 2.7]{Kha18} is nothing else than \cite[Lemma
3.2]{DehKoz12}. Of course, no proof is given for this ``stronger'' version of
\cite[Lemma 3.2]{DehKoz12} in \cite{BejPos19}. Since any $L^{p}$ space with
{$p\in[1,\infty)$} is a modular space, Example \ref{ex1} demonstrates that
\cite[Lemma 1]{BejPos19} is \textbf{false}.

Slight reformulations of \cite[Lemma 1]{BejPos19} can be found in two other
papers by Bejenaru and Postolache; more precisely, \cite[Lemma 2.6]%
{BejPos20b}, attributed to \cite{KhaKoz15} and \cite{DehKoz12} via
\cite{AbdKha17} (although there is no result equivalent to \cite[Lemma
1]{BejPos19} in \cite{AbdKha17}), and \cite[Lemma 1]{BejPos20c}, attributed to
\cite{KhaKoz15}.

\medskip

The list of works with incorrect variations of Schu's lemma
mentioned in this article is by no means complete.

We would like to emphasize that we have only partially verified the
correctness of the results in the aforementioned works that rely on
various variants of Schu's lemma,  such as \cite[Lemma
2.2]{KhAhAb23}, \cite[Lemma 2.4]{KiKiTa02}, \cite[Theorem
2.3.13]{AgOrSa09}, \cite[Lemma 2.2]{ZhoCui15} or \cite[Lemma
1]{BejPos19}.  While some of these works directly apply an incorrect
variant of Schu's lemma in their proofs, others,  though they
strangely cite a false variant,  do verify that the third condition
is met by ``lim''.

\end{document}